\documentclass[fleqn,12pt]{article}
\usepackage{psfig}
\usepackage{amsmath,amstext,amsfonts,amsthm}
\usepackage[mathscr]{eucal}
\newcommand{\nc}{\newcommand}
       
\nc{\tbt}[4]{
  \left[ \begin{array}{cc}
       #1 & #2 \\ #3 & #4
       \end{array} \right] }
\nc{\tbo}[2]{
  \left[ \begin{array}{c}
       #1 \\ #2
       \end{array} \right] }
\nc{\obt}[2]{
  \left[ \begin{array}{cc}
       #1 & #2
       \end{array} \right] }
\nc{\thbo}[3]{
  \left[ \begin{array}{c}
       #1 \\ #2 \\ #3
       \end{array} \right] }
\nc{\thbth}[9]{
  \left[ \begin{array}{ccc}
       #1 & #2 & #3 \\
       #4 & #5 & #6 \\
       #7 & #8 & #9
 \end{array} \right] }

\begin{document}
 
\title{{\bf   
Neural Stabilization/Excitation Control of a High-Order Power System by Adaptive Feedback
Linearization}}

\author{
Kingsley Fregene  \\
Dept. of Electrical and Computer Engineering \\   
University of Waterloo \\
Waterloo, Ontario, Canada N2L 3G1 \\
{\em kocfrege@kingcong.uwaterloo.ca}
\and
Diane Kennedy \\
\thanks{This work was supported by the Natural Sciences and
Engineering Research Council of
Canada via a research grant. }    
Diane Kennedy \\
Dept. of Electrical and Computer Engineering \\
Ryerson Polytechnic University \\
Toronto, Ontario, Canada  M5B 2K3 \\
dkennedy@ee.ryerson.ca
}

\date{}
       
\bibliographystyle{plain}
\maketitle

\begin{abstract}
This paper discusses the systematic design of an adaptive feedback
linearizing neurocontroller 
for a high-order model of the synchronous machine/infinite bus power
system.
The power system is first modelled as an input-output nonlinear
discrete-time system approximated by two
neural networks. The approach allows a simple linear 
pole-placement controller (which is itself not a neural network) to be
designed. The control law is specified
such that the controller adaptively calculates an appropriate feedback
linearizing control law at each sampling instant
by
utilizing
plant parameter estimates provided by the neural system model.
The control system also adapts itself on-line.
This
avoids the
requirement for exact knowledge of the power system dynamics and full
state measurement as
well as other difficulties
associated with  implementing analytical input-output feedback linearizing
control for a complex power system
model.  Furthermore, a departure is made from the `ad hoc' manner in which many neural
controllers have been designed for power systems; the
approach used here has foundations in control theoretic concepts of  adaptive feedback
linearization and pole-placement control design.
 
Simulation results demonstrate the  performance of this controller for a
representative example of a
single-machine/infinite bus power system configuration under various
operational conditions. 

\end{abstract}
%\vspace{10mm}
{\bf Keywords}: Synchronous machines, Neurocontrollers, nonlinear control, Synchronous generator
excitation, Adaptive control, Synchronous generator stability, Power system control, Neural networks, Power
system identification,

\section{Introduction}

The differential geometric technique of state feedback linearization has 
been explored over the past two decades for the purpose of power system control (see, for instance, \cite{diacc},   
\cite{marino}, \cite{mak}, \cite{chapetal}). The main idea is to perform a co-ordinate transformation of the  
nonlinear state space system equations and define a new control input so  that in the new co-ordinates, the nonlinearities in the
plant are either 
wholly or partially masked. This formulation may result in linearizations which are valid for large 
practical operating points of the system, as opposed to a local Jacobian linearization about an operating point.

Feedback linearization has been applied to the control of power systems represented by the state-space single 
machine/infinite bus model in two main ways. The \emph{input-state} feedback linearization approach is formulated 
such that the system state becomes a linear function of a new control input and a new state, while the output (for 
our purposes, the terminal voltage) is still a nonlinear function, e.g. see  \cite{marino}. 
Although this works well for designing stabilizing 
controllers, the nonlinearities in the output map make it very difficult to achieve good tracking of reference 
voltage signals. An alternative to this approach is the \emph{input-output} feedback linearization in which
the terminal voltage 
becomes a linear function of a new control input, e.g. see \cite{diacc}, \cite{dithesis}, \cite{mak}. Due to the mathematical 
complexity of the nonlinear 
state-space model representing the high-order single machine/infinite bus power system (a $7th$ order model), a simplified $3rd$ order
model is 
used in \cite{diacc}, \cite{dithesis}. 
An excitation controller/stabilizer is then designed and tested using the  high order model of the plant. The effects of the 
unmodelled system dynamics and possible shifts in parameters of the plant
raise some questions on the effectiveness of the controller.
An inherent
drawback of the feedback linearization approach is the non-robustness  due to the need for exact knowledge, both in
terms of the structure
and  parameters, of the 
model of the power system. Adaptive control has been used to compensate for parameter variation in a framework
that allows
the controller to learn the nonlinearities on-line \cite{isidori2}. However, this assumes that the nonlinearities can be
parameterized linearly in some unknown parameters \cite{isidori2}, \cite{koko}. Some
power system stabilizers based on adaptive control design
techniques have been shown to increase the operating range over
which they can provide good control as well as exhibit robustness
to parameterized system disturbances \cite{pierre_87}.

This paper presents a design which addresses these problems by working directly with the more 
accurate $7th$ order model of the power system. The nonlinear dynamical system is modelled using multilayer 
neural networks made of sigmoid type nonlinearities and adjustable weights which are nonlinearly parameterized. The flexibility provided by
the nonlinear parameterization allows for a 
more realistic representation of the underlying nonlinear power system model. Power 
system controller designs using neural networks of various configurations exist in the literature and are reported in 
works like \cite{vidya}. Neural adaptive control of feedback linearizable nonlinear systems was first proposed in 
\cite{chen1} and extensively analyzed in \cite{chen2}. In these references, the networks were trained using the  backpropagation 
algorithm 
\cite{rumcle}; the nonlinear function approximation capability of neural networks (see \cite{funahashi}) was exploited to 
model the plant in order to determine an appropriate feedback linearizing control. However, the approach used in \cite{chen1},
\cite{chen2} provides no clear way of specifying closed loop dynamics.  
In this work the neural network identifier is trained in batch form and adapted using the 
Levenberg-Marquardt optimization (as in \cite{magnus1}) rather than using standard backpropagation. Also, the whole control problem is
formulated such that tracking dynamics can be user-specified by approprately assigning closed loop system poles. Weight adaptation
based on nonlinear optimization is attractive here because it significantly corrects some of the defects of standard backpropagation, 
most importantly slow rate of convergence. In this framework, 
the trained neural networks give an on-line estimate of the parameters of the power system at every sampling instant which are 
then used to calculate the feedback linearizing control law. The controller is designed to
provide 
reference tracking (with bounded state), damping of power angle oscillations and tolerate variations in 
the plant parameters. It may also be adapted on-line.

This paper is organized as follows: Section 2 presents the state-space synchronous machine/infinte bus power system model. The system 
relative degree and minimum phase properties, important for feedback linearization to be successfully applied, are also discussed.
The class of systems to which 
this model belongs and which motivates the choice of a neural identifier is determined. 
Section 3 focuses on the neural modelling of the power system while Section 4 presents a systematic design of the proposed excitation
controller/power system stabilizer. 
Finally, simulation results are presented and discussed in  
Section 5, followed by some concluding remarks.     

\section{System and Identification Model} 
\begin{figure}[h]
\begin{center}
\input{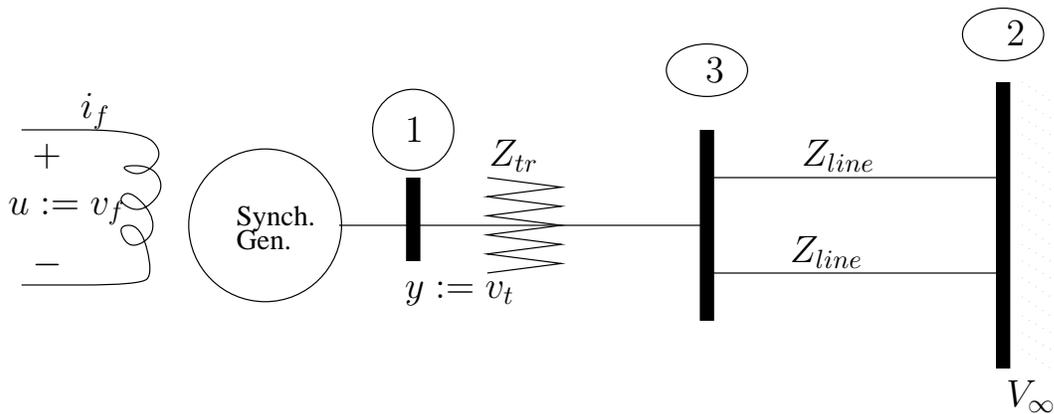}
\caption{Synchronous Machine Connected to an Infinite Bus}
\label{genlin}
\end{center}
\end{figure}

The model used is developed in \cite{kundur} (p. 54ff). It
yields a Park's 7th order nonlinear 
time-invariant, state-space description of the synchronous generator/infinite bus system. This is connected through a balanced pair of transmission 
lines to a power system modelled simply as one of infinitely large ability to supply and absorb power at a set voltage level, 
i.e. an \emph{infinite bus} system. The system is depicted in Fig. \ref{genlin}. 

State variables for this generator/infinite bus set-up are selected as follows:

\begin{itemize} 
\item $\delta (x_{1})$ - power angle
\item $\omega (x_{2})$ - power angle derivative
\item $\lambda_{d} (x_{3})$ - d-axis flux linkage
\item $\lambda_{q} (x_{4})$ - q-axis flux linkage 
\item $\lambda_{f} (x_{5})$ - field flux linkage
\item $\lambda_{k_{d}} (x_{6})$ -  damper d-axis flux linkage
\item $\lambda_{k_{q}} (x_{7})$ -  damper q-axis flux linkage 
\end{itemize}
The system model is given by
\begin{eqnarray}
\dot \delta & = & \omega    \nonumber \\
\dot \omega & = & \frac{\omega_{b}}{2H}[P_{m}-P_{e}(\delta,\lambda_{d},\lambda_{q}) -D \omega] \label{begin_nonlinear} \\
\dot \lambda & = & \omega_{b} \left( [RL^{-1} + Z] \lambda + [v_{d} \; v_{q} \; v_{f} \; 0 \; 0]^{T} \right) \nonumber 
\end{eqnarray}
where $\lambda := [\lambda_{d} \; \lambda_{q} \; \lambda_{f} \; \lambda_{k_{d}} \; \lambda_{k_{q}}]^{T}$, are the per 
unit flux linkages, R := $diag \; [r_{s} \; r_{s} \; -r_{f} \; -r_{k_{d}} \; -r_{k_{q}}],$ are per unit resistances. $H$ is the inertia constant, $D$
encapsulates the damping in the system,

\begin{eqnarray*}
L :=
\left[
\begin{array}{ccccc} -L_{d}&0&L_{ad}&L_{ad}&0 \\ 0&-L_{q}&0&0&L_{aq} \\ -L_{ad}&0&L_{f}&L_{fk_{d}}&0 \\
                     -L_{ad}&0&L_{fk_{d}}&L_{k_{d}}&0 \\ 0&-L_{aq}&0&0&L_{k_{q}} \end{array} \right] 
\end{eqnarray*}
are inductances, and
\begin{eqnarray*}
Z :=
\left[
\begin{array}{rcccc} 0&1&0&0&0 \\ -1&0&0&0&0 \\ 0&0&0&0&0 \\ 0&0&0&0&0 \\ 0&0&0&0&0 \end{array}
\right]_{.}
\end{eqnarray*}
As shown in Fig. \ref{genlin}, the control input ($u$) is the field voltage $v_{f}$ and the output to be
controlled ($y$) is the
generator terminal
voltage $v_{t}$ given by
\begin{eqnarray*}
y=v_{t}=h(\delta,\lambda_{d},\lambda_{q})=\sqrt{v^{2}_{d}+v^{2}_{q}}
\end{eqnarray*} 
where
\begin{eqnarray*}
v_{d} & = & r_{11}I_{d}-x_{11}I_{q} + v_{\infty} [A cos (\pi/2-\delta) + B sin(\pi/2 - \delta)] \\  
v_{q} & = & r_{11}I_{q}-x_{11}I_{d} - v_{\infty} [B cos (\pi/2-\delta) - A sin(\pi/2 - \delta)]_{;}
\end{eqnarray*}
the variables are described in \cite{kundur} and \cite{dithesis}. $A, \; B$ are constant matrices having values depending on physical 
parameters of the power system, 
$v_{\infty}$ is the voltage of the infinite bus, $I_{d,q}$ are $d-$ and $q-$ axis currents, $r_{11}$, $x_{11}$ are components 
of transmission line resistance and reactance (see 
\cite{kundur}). The system model 
(\ref{begin_nonlinear}) belongs to the class of control affine nonlinear systems of the form 
\begin{eqnarray}
\dot{x} & = & f(x)+g(x)\cdot u \label{affine_sys} \\   
y & = & h(x) . \nonumber
\end{eqnarray}
Here  $f$, $g$ are smooth
functions; the state $x \in \mathbb{Re}^{n}$, 
input $u \in \mathbb{Re}$ and output $y \in \mathbb{Re}$. The \emph{relative degree} and \emph{minimum phase} properties of this system 
are examined next.

The relative degree of the power system (\ref{begin_nonlinear}) may be
obtained by differentiating  $y$ (i.e. $v_{t}$)  
until $u$ ($v_{f}$) appears. This is achieved by computing the \emph{Lie
derivatives} (see \cite{isidori}) of the 
system.  
\begin{eqnarray*}
L_{f}h & = & \nabla h \dot{x} \\  
       & = & \frac{dh}{dx_{1}} \dot{x}_{1} + \ldots + \frac{dh}{dx_{5}} \dot{x}_{5} + \ldots + \frac{dh}{dx_{7}} \dot{x}_{7} 
\end{eqnarray*}
Observe that  $v_{f}$ appears explicitly in the expression for $\dot{x}_{5}(\dot{\lambda}_{f})$ in (\ref{begin_nonlinear}). 
Therefore the 
output need only be differentiated once before the control input $v_{f}$ appears, i.e. {\em the relative degree is 1}.
It is easy to show that the system (\ref{begin_nonlinear}) is nominally {\em minimum phase} over the practical operating region by
linearizing the system at various
operating points in the state-space and verifying left half plane
zeros (see Appendix B in \cite{kofthesis}).

From the relative degree and minimum phase properties, the identification model proposed for the system will be the 
single-input single-output relative-degree-one nonlinear discrete-time system
\begin{eqnarray}
y[k+1] & = & f[y(k),\ldots,y(k-n+1),\;u(k-1),\ldots,u(k-m)] \; +  \label{ident_eqn} \\ 
        &   & g[y(k),\ldots,y(k-n+1),\;u(k-1),\ldots,u(k-m)] \cdot u(k) \nonumber
\end{eqnarray}
where $y$ is the output, $u$ is the input, $f$ and $g$ are assumed smooth and the system is
assumed minimum phase. For the $7^{th}$-order power system  with relative degree 1, n = 7, m = 6. 

\section{Modelling the power system} 
\label{sec_modelling}
\begin{figure}[h]
\centering\
\psfig{file=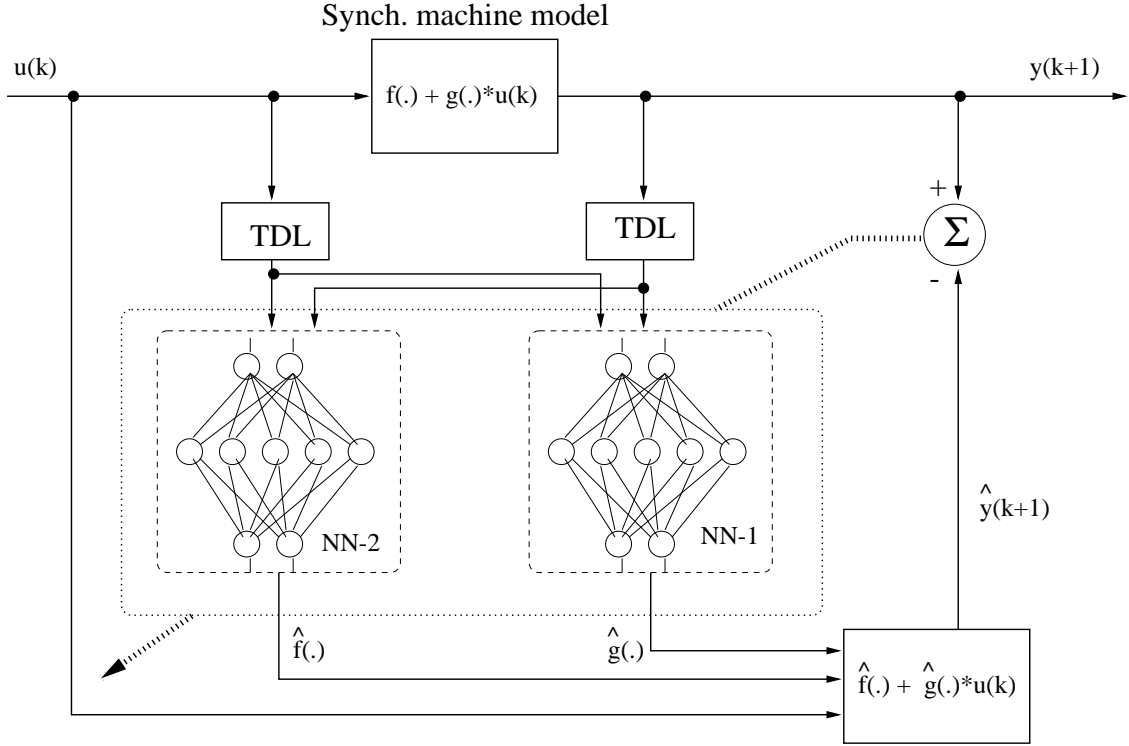,width=15cm}
\caption{Neural system identification of the synchronous machine model}
\label{ident}
\end{figure}

The first stage in the controller design is to identify the synchronous machine/infinite bus dynamical model using the structure 
(\ref{ident_eqn}). The set-up to do so is shown in Fig. \ref{ident}. $TDLs$ are tapped delay lines 
used to implement the regressor structure for both networks.
\emph{NN-1} and \emph{NN-2} are feedforward neural networks whose input vector consists of present output
and delayed inputs and outputs. They are used to predict the next output $y(k+1)$. 
\emph{NN-1} has one hidden layer made up of a total of $p$ neurons with
hyperbolic tangent \emph{tanh} activation function and
an output layer of one neuron with linear activation function. \emph{NN-2}
has $q $ \emph{tanh} hidden
layer neurons and one output neuron. All weights are initialized randomly and 
contain biases.

A training set is generated by exciting the plant with a variety of  input signals
$u$ (of appropriate range) and measuring the output $y$. That is, let $Z^{N}=\{ [u(k), y(k)]|k=1,\ldots,N \}$, represent the training
data pairs and 
N be the total number of data sets. The estimate of the plant output given
by the neural nets is
\begin{equation}
\hat{y}[k+1] = \hat{f}[z(k),{\bf w}] + \hat{g}[z(k),{\bf v}] \cdot u(k)  \label{nn_out}
\end{equation}
where
\begin{eqnarray*}
\hat{f}(z(k),w) & = & \sum^{p}_{i=1}w_{i}\: \tanh \left(\sum^{m+n}_{j=1} \:
w_{ij}z_{j}+b_{i} \right) + b  \\
\hat{g}(z(k),v) & = & \sum^{p}_{i=1}v_{i}\: \tanh \left(\sum^{m+n}_{j=1} \:
v_{ij}z_{j}+a_{i} \right) + a  ;
\end{eqnarray*}
$(w,b)$ and $(v,a)$ are weights and biases respectively for the $\hat{f}(\cdot)$ and $\hat{g}(\cdot)$
networks. Let the set of possible weights be $\hat{\Theta}$ and define the cost
function ${\bf J}$ as the mean square error in predicting the plant
output. That is,
\begin{equation} J_{N}(\Theta,Z^{N}) = 
%\begin{equation}
\frac{1}{2N}\sum^{N}_{k=1} \left[y(k+1)-\hat{y}(k+1|\Theta)\right]^{T}
\left[y(k+1)-\hat{y}(k+1|\Theta)\right].    \label{cost_1}
\end{equation}
The optimal weights are determined by solving
\[ \hat{\Theta} = arg \; \min_{\Theta} \; J_{N} (\Theta,Z^{N}) \]
Let the weight update rule toward achieving the optimal choice be given by
\[ \Theta_{k+1} = \Theta_{k} + \rho_{k} S_{k} \]
where $\Theta_{k}$ is the current iterate of the set of weights, $\rho_{k}$ is the step size and $S_{k}$
is the search direction.

For identification purposes, the weight update is done in batch form using 
the Levenberg-Marquardt optimization as described in \cite{magnus1}. When 
the cost function is minimized to an acceptable level, weight adaptation is
discontinued and the neural model is validated on the training data set and
subsequently cross-validated on other sets of data not in the original training set.
Simulation results are provided which show the effectiveness of the proposed identification technique.

\section{Adaptive linearizing control design}
\label{adapt_lin_contr}
Recall that the power system  is modelled by the SISO discrete time system:
\begin{eqnarray}
y[k+1] & = & f[y(k),\ldots,y(k-n+1),\;u(k-1),\ldots,u(k-m)] \; +  \label{ident_eqn_cont} \\
       &   &  g[y(k),\ldots,y(k-n+1),\;u(k-1),\ldots,u(k-m)] \cdot u(k) \nonumber
\end{eqnarray}
and estimated by the neural model:
\begin{eqnarray}
\hat{y}[k+1] & = & \hat{f}[y(k),\ldots,y(k-n+1),\;u(k-1),\ldots,u(k-m)] \; +  \label{neural_ident_eqn} \\
             &   & \hat{g}[y(k),\ldots,y(k-n+1),\;u(k-1),\ldots,u(k-m)] \cdot u(k) \nonumber
\end{eqnarray}
A simple pole-placement controller is proposed to allow the tracking dynamics to be specified by
appropriately placing closed-loop poles at stable locations.
Let the reference signal  to be tracked be given by $r(k)$.

To facilitate the discussion, a state-space form of the system (\ref{ident_eqn_cont})
is constructed by choosing state variables as follows:
\begin{eqnarray*}
z_{11}(k) & = & y(k-n+1) \\
          & \vdots   &  \\
z_{1n}(k) & = & y(k) \\
z_{21}(k) & = & u(k-m) \\
          & \vdots   & \\
z_{2m}(k) & = & u(k-1) 
\end{eqnarray*}
and we subsequently define
\begin{eqnarray*}
{\bf z} & = & \left[ \begin{array}{c} {\bf z}_{1} \\ {\bf z}_{2} \end{array} \right] \;  = \; \left[ \begin{array}{c} z_{11}\\ \vdots \\
z_{1n} \\ -   
\\ z_{21} \\ \vdots \\ z_{2m} \end{array} \right]_{.}
\end{eqnarray*} 
This leads to the state-space model:
\begin{eqnarray}
z_{11}(k+1) & = & z_{12}(k) \nonumber \\
            & \vdots &  \nonumber \\
z_{1n}(k+1) & = & f[{\bf z}(k)] + g[{\bf z}(k)] \cdot u(k)=y[k+1]
\label{ssout_model}\\
z_{21}(k+1) & = & z_{22}(k) \nonumber \\
            & \vdots & \nonumber \\
z_{2m}(k+1) & = & u(k) \nonumber \\
y(k) & = & z_{1n}(k) \nonumber
\end{eqnarray}
Now we define a new control input $\tilde{u}$ and set
\begin{eqnarray} 
u(k) & = & \frac{1}{g[{\bf z}(k)]}\left[-f[{\bf z}(k)]+\tilde{u}(k)\right];
\label{fblc_u}
\end{eqnarray}
observe that this results in input-output closed loop dynamics of
\begin{eqnarray*}
y[k+1] & = & \tilde{u}(k)
\end{eqnarray*}
Now we address the primary control objective of tracking a constant reference voltage $r=v_{t}$.
This motivates the choice of the feedback control $\tilde{u}(k)$ as
\begin{eqnarray} 
\tilde{u}(k)=K_{1}r(k) - \left[ C_{p-1}y(k)+C_{p-2}y(k-1)+ \cdots + C_{0}y(k-p+1) \right]
\label{ytilde}
\end{eqnarray}
with $p \in \mathbb{Z}^{+}$ and $C_{i} \in \mathbb{R}$. Hence,
\begin{eqnarray*}
y[k+1] +  C_{p-1}y(k)+C_{p-2}y(k-1)+ \cdots + C_{0}y(k-p+1) & = & K_{1} r(k).
\end{eqnarray*}
Therefore the transfer function between the reference and $y$ is given by
\begin{eqnarray}
\frac{V_{t}[z]}{V_{ref}[z]} &=& \frac{Y[z]}{R[z]}  =  \frac{K_{1}}{z^{p}+C_{p-1}z^{p-1}+\cdots+C_{1}z+C_{0}}
\; =: \; \frac{K_{1}}{Q(z)}.
\label{transfer_function}
\end{eqnarray}
Closed loop poles for the system can be assigned by selecting constants $C_{i}$'s so that the zeros of the characteristic polynomial $Q(z)$
are all in
the unit circle.

Some assumptions made on the system and a discussion of the stability of this kind of control scheme can be found in \cite{kofthesis}.
Specifically, it is shown therein that the approach used here is a more general form  of that extensively analyzed in \cite{chen2}.

Now we rewrite the plant in the input-output form as
\begin{equation}
\label{newplant}
y(k+1) \; = \; f[{\bf z}(k)] + g[{\bf z}(k)] \cdot u(k)
\end{equation}
The plant (\ref{newplant}) is
thus modelled by the neural networks as:
\begin{equation}
\label{newplant_neural}
\hat{y}(k+1) \; = \; \hat{f}[{\bf z}(k),{\bf w}] + \hat{g}[{\bf z}(k),{\bf v}] \cdot u(k)
\end{equation}
For the three-layer neural networks used in the identification stage, with total
hidden
neurons $p$ and $q$ respectively and weights/biases $(w_{i}, \; \hat{w}_{i} \in {\bf w})$ and $(v_{i}, \; \hat{v}_{i} \in {\bf v})$ , the
functions $\hat{f}(\cdot,\cdot)$ and $\hat{g}(\cdot,\cdot)$  can be expressed as
\begin{equation}
\hat{f}[{\bf z}(k),{\bf w}]=\sum^{p}_{i=1}w_{i}\: \tanh \left(\sum^{m+n}_{j=1} \:
w_{ij}x_{j}+\hat{w}_{i} \right)
\end{equation}
and
\begin{equation}
\hat{g}[{\bf z}(k),{\bf v}]=\sum^{p}_{i=1}v_{i}\: \tanh \left(\sum^{m+n}_{j=1} \:
v_{ij}x_{j}+\hat{v}_{i} \right)
\end{equation}
Funahashi's theorem \cite{funahashi} says that we can approximate $f$ and $g$ as close as we like  by a suitably
parameterized
neural network. It is not presently known how to choose the number of neurons given a desired error tolerance. However, the typical
approach
is to choose a configuration and try to obtain good weights/biases by trying to solve a nonlinear optimization problem. If the
approximation
is poor, the configuration is adjusted and repeated. This will be accomplished via an ``off-line'' training phase using the
Levenberg-Marquardt algorithm (as discussed in Section \ref{sec_modelling}), to obtain initial estimates of $\hat{f}$ and $\hat{g}$.

Now, we focus on the on-line phase, in which the controller is applied to the system. An additional constraint to consider is the need to
ensure that $\hat{g}$ remains nonzero due to the fact that the controller incorporates a division by $\hat{g}$ as shown in
(\ref{neur_cont_law}).

\begin{figure}[h]
\centering\
\psfig{file=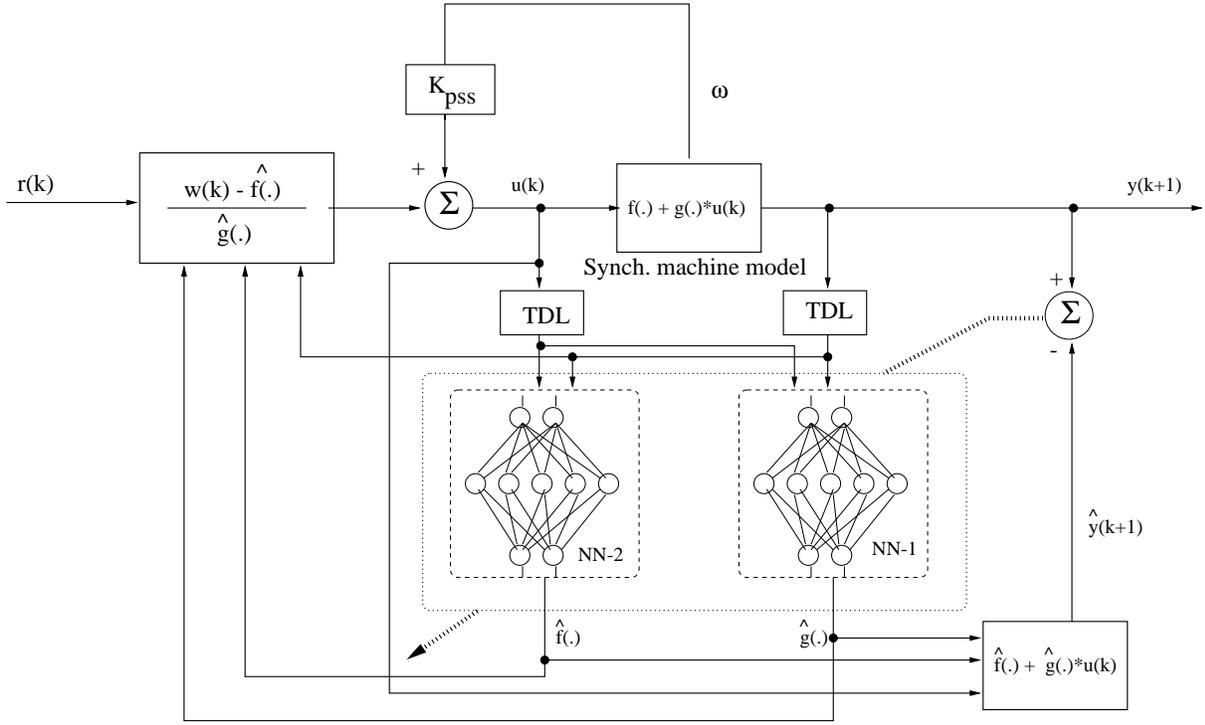,width=16cm}
\caption{Adaptive linearizing control/power system stabilizer set-up}
\label{cont_pss}
\end{figure}

Let ${\bf w}(k)$ and ${\bf v}(k)$ denote the estimates of ${\bf w}$ and ${\bf v}$ at sampling instant  $k$. The
control input 
$u(k)$ is then given as: \\ \\
{\em \bf{Control Law:}}
\begin{equation}
\vspace{10mm}
\label{neur_cont_law}
u(k) \; = \; \frac{- \hat{f}[{\bf z}(k),{\bf w}(k)] + \tilde{u}(k)}{\hat{g}[{\bf z}(k),{\bf v}(k)]}
\end{equation} \\
with $\tilde{u}(k)$ as previously defined in (\ref{ytilde}). The control $u(k)$ is applied to both the plant and the neural network model;
the network weights are updated using the
error between the plant output and the output predicted by the neural model. At each sampling instant $k$, the estimated
plant output is given by
\begin{equation}
\label{estimate_op}
\hat{y}^{*}(k+1) \; = \; \hat{f}[{\bf z}(k),{\bf w}(k)] \; + \; \hat{g}[{\bf z}(k),{\bf v}(k)] \cdot u(k)
\end{equation}
The error associated with this estimate $\{ e^{*}(k+1) \}$ is defined as
\begin{equation}
\label{error_star}
e^{*}(k+1) \; = \; \hat{y}^{*}(k+1) - y(k+1)
\end{equation}
This error (together with a {\em deadzone} adaptation) will be used in the weight update rule described next. 
A deadzone is an absolute lower limit on the adaptation error below which no adaptation should be done. That is to say, before the
weights are adapted at each sampling instant, the plant output error $e^{*}(k+1)$ is compared with the dead zone radius
$d_{0}$. If $\left|e^{*}(k+1)\right| \leq d_{0}$, no adaptation is done, otherwise the adaptation is done using an error adaptation
rule based on the deadzone $D(e)$ given by
\begin{eqnarray}
\label{d_zone_eqn}
D(e) =
\left\{
\begin{array}{lcl}
0 & if & \left| e \right| \leq d_{0} \\
e-d_{0} & if & e>d_{0} \\
e+d_{0} & if & e<-d_{0} \\
\end{array}
\right.
\end{eqnarray}
which ensures that the constraint of nonzero $\hat{g}$ holds.
In order to implement the weight update rule, first define the set of network weights $\Theta$ as
\[  \Theta \; = \; \left[ \begin{array}{cc} {\bf w} \\ {\bf v} \end{array}
\right] \] \\
and the Jacobian to be computed at each time step by
\begin{eqnarray}
\label{online_cost}
{\bf J}_{k} & = & \left[ \left. \frac{\partial \hat{y}^{*}(k+1)}{\partial \Theta} \right|_{\Theta(k)} \right]^{'} \nonumber \\
                   \nonumber \\
            & = & \left[ \begin{array}{c} \left( \left. \frac{\partial \hat{f}[{\bf z}(k), {\bf w}]}{\partial {\bf w}}
\right|_{{\bf w(k)}} \right)^{'} \\
\left( \left. \frac{\partial \hat{g}[{\bf z}(k), {\bf v}]}{\partial {\bf v}}
\right|_{{\bf v(k)}} \right)^{'} \cdot u(k)
\end{array}
\right]_{.}
\end{eqnarray}
This Jacobian can be calculated using the backpropagation algorithm \cite{rumcle} at the end of each
iteration, since the  
variable
$\hat{y}^{*}(k+1)$ is the output of the neural networks, and is available for use \cite{chen2}.\\
\vspace{6mm}
{\em {\bf Adaptation Rule:}} \\
\begin{equation}
\label{online_update_rule}
\Theta(k+1) \; = \; \Theta (k) - \frac{1}{(1+{\bf J}'_{k}{\bf J}_{k})}D\left[e^{*}(k+1)\right]{\bf J}_{k}  
%\vspace{15mm}
\end{equation}
%\vspace{10mm}
%-------commented out on 10/18/99---------
%The implementation  of the on-line control scheme can be summarized  as  follows:
%\begin{enumerate}
%\item $k=0$   
%\item read initial weights for the $\hat{f}$ and $\hat{g}$ networks
%\item sample reference signal, $r(k)$ \label{get_ref}
%\item compute $\tilde{u}(k)$
%\item compute control signal $u(k)$ \label{get_u}
%\item apply $u(k)$ to both the actual plant and the neural network
%\item compute output predicted by the neural network, $\hat{y}^{*}[k+1]$
%\item sample output from the actual power system model, $y[k+1]$
%\item calculate the output error $\epsilon^{*} (k+1) = \hat{y}^{*}[k+1] - y[k+1] $
%\item {\bf if} $|e^{*} (k+1)| \leq d_{0}$ \\ \hspace{10mm} {\bf then} \\ \hspace{20mm} $k
%\longleftarrow k+1$ \\ \hspace{20mm} go to
%Step \ref{get_ref} \\{\bf else} \\
%\hspace{20mm} $\left[
%e^{*} (k+1) \longleftarrow
%      \left( | e^{*}(k+1)| - d_{0} \right)  \right]$ , \\ \hspace{20mm} {\em update weights} \\
%\hspace{20mm}
%      $k \longleftarrow k+1$  \\  {\bf then} \\ \hspace{20mm}
%      go to Step \ref{get_ref}
%\item continue till last reference sample is received
%\end{enumerate}
%\vspace{5mm}
{\em {\bf Power system stabilizer:}}
The {\em power system stabilizer} adds damping to the generator rotor angle oscillations by modulating its excitation using
auxiliary stabilizing signal(s). In order to accomplish this, the stabilizer should be desgined to produce a component of electrical
torque in phase with the rotor speed deviations. For the machine model being used, the signal of choice for modulating the
excitation is the rotor angle
derivative $\dot{\delta}$. This requires that the linearizing control law (\ref{fblc_u}) be modified to introduce a torque
component     
which damps out the oscillations of the rotor angle. The new control law is given by
\begin{equation}
u(k)=\frac{\tilde{u}(k)-\hat{f}(\cdot)}{\hat{g}(\cdot)} + K_{pss}\cdot \dot{\delta} (k)   \label{fblc_u_pss}
\end{equation}
where $K_{pss}$ is a suitable weighting of the rotor angle derivative.
The gain $K_{pss}$ for the neural controller is selected using a root locus criterion
extensively discussed in \cite{dithesis}.
The proposed excitation
controller/power system stabilizer control set-up is
depicted in Fig. \ref{cont_pss}.

\section{Simulation Results and Discussions} 
\begin{figure}[h]
\centering\
\psfig{file=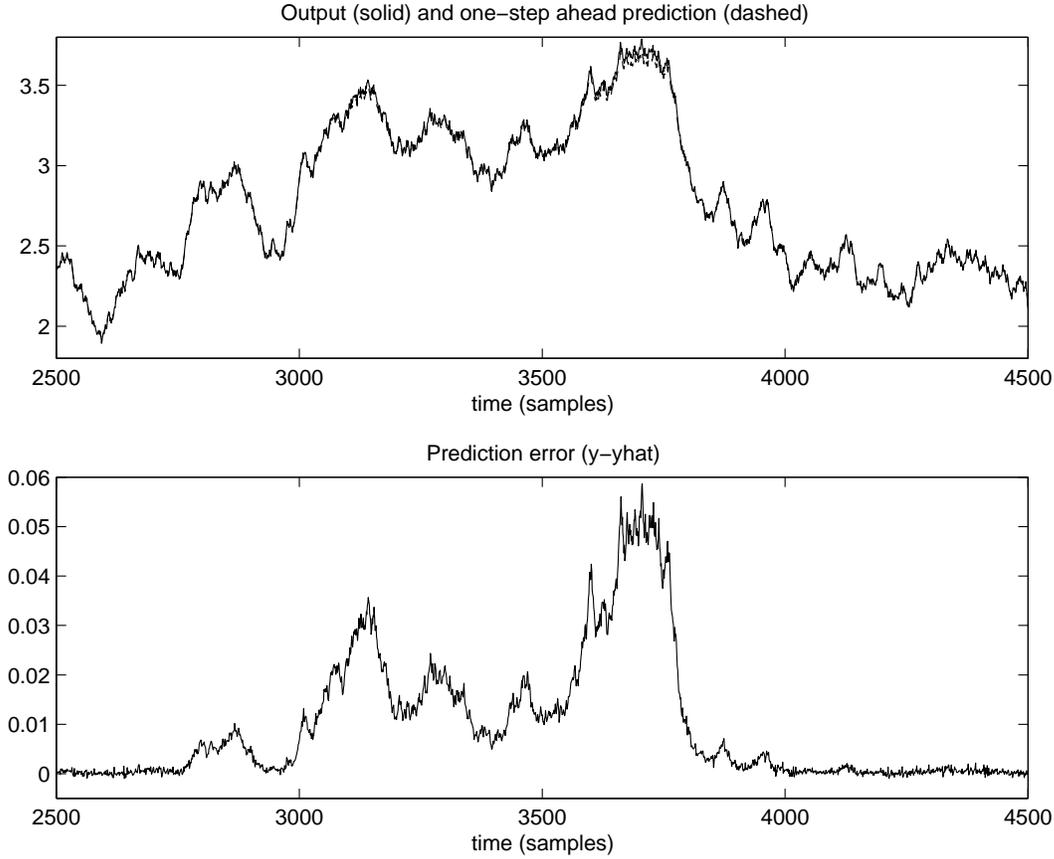,width=14cm}
\caption{Validating the neural identification scheme on test data}
\label{valid_test}
\end{figure}

In order to accomplish the system identification part, the power system is modelled in SIMULINK and sampled at 
2 ms to ensure the fast electrical dynamics of the plant are
adequately captured. Next, 
random input signals $v_{f} \in [-0.1 \; 0.1]$ are used to excite the plant and the terminal voltage $v_{t}$ 
is measured for each case. 10,000 such data pairs $(v_{f},\; v_{t})$ are collected and  divided 
into two parts - one for training 
the neural nets and the other for cross-validating the resulting neural model. The neural nets have the following structure
\begin{itemize}
\item $\hat{f}(\cdot)$ and $\hat{g}(\cdot)$ both have 5 hidden hyperbolic tangent neurons with bias and one linear output neuron
\item The weight matrices for both $\hat{f}(\cdot)$ and $\hat{g}(\cdot)$ have dimension $5 \; \times \; 14$ in the hidden 
      layer and $1 \; \times \; 6$ in the output layer.
\item The input to each network consists of current and delayed $v_{f}$ and $v_{t}$ up to the most delayed signal 
      in each case (13 in all) with the bias term accounting for the last column of the weight matrix  
\end{itemize}
The optimization routine used for offline weight adaptation is allowed to proceed for a total of 150 iterations at the end of which
the cost function  is minimized to the order of $10^{-6}$.
Fig. \ref{valid_test} (top) shows the cross-validation of the model 
on a portion of test data not included in the original training set while Fig. \ref{valid_test} (bottom)
shows the prediction
error on 
this data (from this, the deadzone radius $d_{0}$  is selected as $0.01$). Figures showing the validation of the neural model on the
training
data are not included since they 
do not, in general, give a reliable indication of how well the neural model performs.   
Observe from the error plot that the error is essentially  close to zero
for most of the samples but is more pronounced around the
range 3000 - 4000; even in this range, the errors
are not unduly large. The maximum value of $\epsilon^{*}(k+1)$ in this range is $0.06$ resulting in a small relative error of
\[
\frac{ \left| \epsilon^{*} \right|_{\infty}}{ \left| y \right|_{\infty}} \;  = \; \frac{0.06}{2.25} \times 100 \% \; = \; 2.7 \%
\]
Observe also that this is the maximum error and it only occurs for a short sampling range.
We argue that adapting the weights of the network on-line will serve to compensate for the effects of this
error on controller performance.

We begin with voltage tracking only and compare the neural excitation controller with both a conventional ST1A high gain excitation controller (see \cite{kundur} p. 365) and
a nonlinear analytic feedback linearizing excitation controller proposed in \cite{dithesis}. The gain and time constant  of the ST1A
excitation controller are
$K_{e} =
200$ and $ T_{e} = 0$, so that the control input is given by
\[ v_{f} \; = \; 200 \frac{r_{f}}{x_{ad}} (v_{ref} - v_{t}). \]
All quantities are in per unit and
\[ \frac{r_{f}}{x_{ad}} = 3.9056 \times 10^{-4} \]
so that the control law is given by
\begin{equation}  
\label{exciter_u}
v_{f} \; = \; 0.0781 (v_{ref} - v_{t})
\end{equation}
The ST1A exciter is based on a linearization at the operating point $v_{ref}=1.1392$. 

The analytic excitation controller is designed in \cite{dithesis}  and has a control of the form
\begin{equation}
v \; = \; -K(v_{t} - v_{ref}), \hspace{10mm} K > 0
\end{equation}
where $v$ is a virtual control resulting from performing input-output feedback linearization on a reduced order model of
the power system. The closed loop transfer function is given by
\[ \frac{v_{t}(s)}{v_{ref}(s)} \; = \; \frac{K}{s+K} \]
so that the terminal voltage tracks the reference exponentially at a rate depending on the controller gain $K$. For the
simulations done here
\[ K \; = \; 25 \]
in order to provide a  transient response similar to that of the ST1A controller.

Recall that the neural excitation control law is given by
\[ v_{f}(k) \; = \; \frac{-\hat{f}(\cdot) + \tilde{u}(k)}{\hat{g}(\cdot)} \]
where
\begin{eqnarray}
\tilde{u}(k)=K_{1} r(k) -
              \left[ C_{p-1} \; C_{p-2} \; \ldots \; C_{0} \right] \cdot \left[ \begin{array}{c} y_{k} \\ y_{k-1}\\ \vdots
\\
              y_{(k-p+1)} \end{array} \right]_{.}
\label{more_ytilde}
\end{eqnarray}
and $r(k)$ is a constant reference voltage.
The constants $C_{i}$ are obtained from the characteristic polynomial of the desired closed-loop system given by
\[ Q(z) \; = \; z^{p}+C_{p-1}z^{p-1}+,\ldots,+C_{1}z+C_{0} \]
and $K_{1}$ is chosen to be $1 + \sum_{i=0}^{p-1} C_{i}$, to yield unity DC gain.
The transfer function is then given by
\[ \frac{V_{t}[z]}{V_{ref}[z]}= \frac{Y[z]}{R[z]} = \frac{K_{1}}{z^{p}+C_{p-1}z^{p-1}+\cdots+C_{0}}  \]
It remains to select $p$ appropriately so that `nice' control action is obtained. 

\newpage
{\bf Some special cases}
\begin{itemize}
 \item {\em p=0}: \hspace{2mm}
Observe that setting $p=0, \; K_{1}=1$ causes (\ref{more_ytilde}) to reduce to
\[ \tilde{u}(k) \; = \; r(k) \]
and this is how the reference signal is obtained in the work of Chen and Khalil \cite{chen2}. This approach
gives a control structure which lends itself to stability analysis and allows parameter convergence to be
guaranteed.
However, the dynamics of the system can not be prescribed (for instance by trying to place poles of the
closed loop system). Furthermore, using this scheme for the excitation control/pss application shows
a harsh demand for tracking and very high level of control action.
The control set-up in this case is shown in Fig. \ref{fig_chen1} while simulation result is depicted
in Fig. \ref{case_p0}.
\begin{figure}[h]
\centering\
\psfig{file=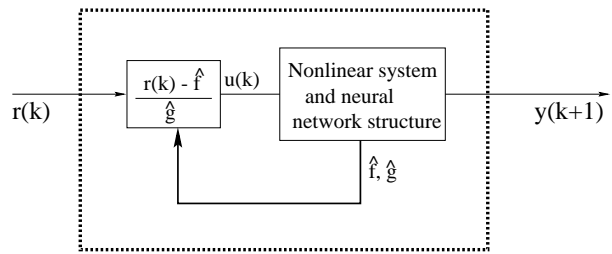,width=8cm}
\caption{Excitation Control setup for Case $p=0$}
\label{fig_chen1}
\end{figure}

\begin{figure}[h]
\centering\
\psfig{file=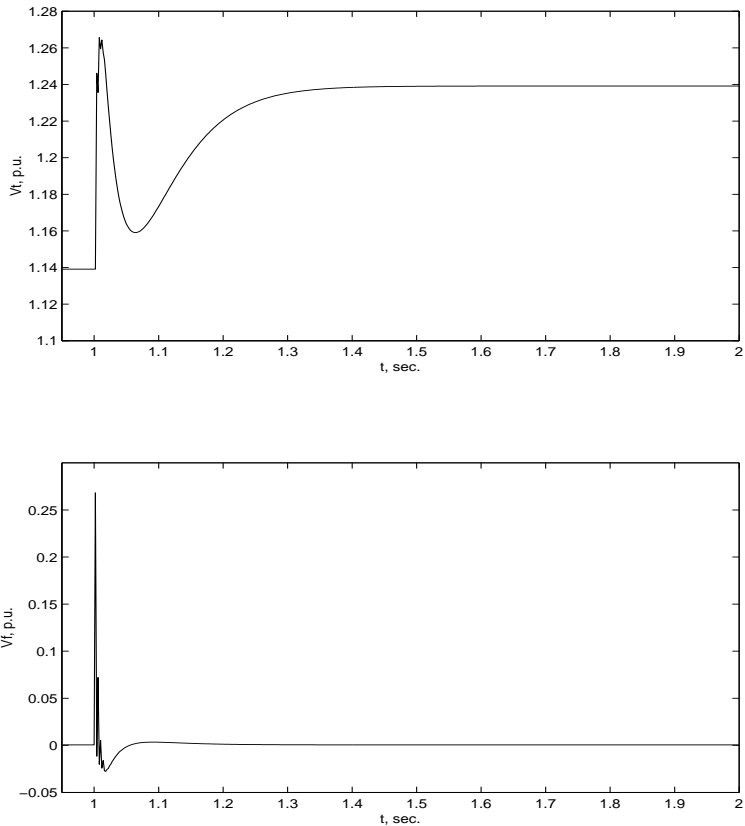,width=10cm,height=11cm}
\caption{Excitation Controller for Case $p=0$}
\label{case_p0}
\end{figure}

\item {\em Case $p=1$}: \hspace{2mm} Setting $p=1$ causes
(\ref{more_ytilde}) to
reduce to
\begin{eqnarray} \tilde{u}(k) \; =  \; (1 + C_{0}) r(k) - C_{0} \cdot y_{k} \label{casep1} \end{eqnarray}
so that the transfer function of  the resulting closed loop system is
\[ \frac{Y[z]}{R[z]} = \frac{1+C_{0}}{z+C_{0}} \]
with a first order characteristic polynomial. This should give a `nice' first order tracking behaviour and
constrain the level of control
effort.
To check the veracity of this reasoning, the controller/pss is tested with $\tilde{u}$ defined as in
(\ref{casep1}) and $C_{0}$ selected to give similar speed of response as the ST1A conventional controller
used for comparison. The control setup in
this case is
shown in Fig. \ref{fig_diane1} while the response is depicted in Fig.
\ref{case_p1}, where it is compared to the analytic controller design based on a reduced order model of the
machine reported in
the work of Kennedy, {\em et. al.}  \cite{dithesis}.
\begin{figure}[h]
\centering\
\psfig{file=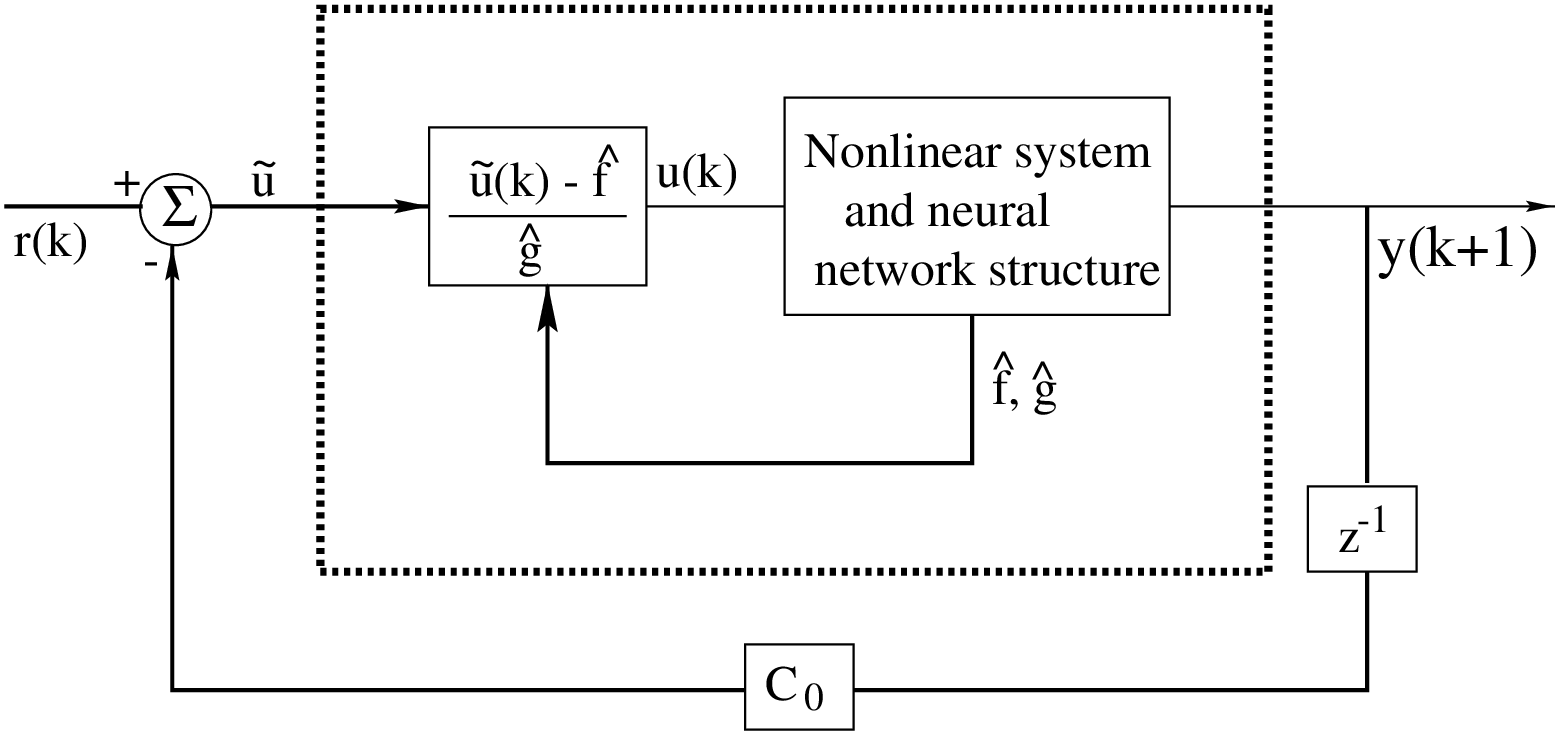,width=8cm}
\caption{Excitation Control setup for Case $p=1$}
\label{fig_diane1}
\end{figure}
Observe that the neural controller is no better than the analytic controller in this case.
Furthermore, it is apparent that the response is not that of a first order linear system. 
\begin{figure}[h]
\centering\
\psfig{file=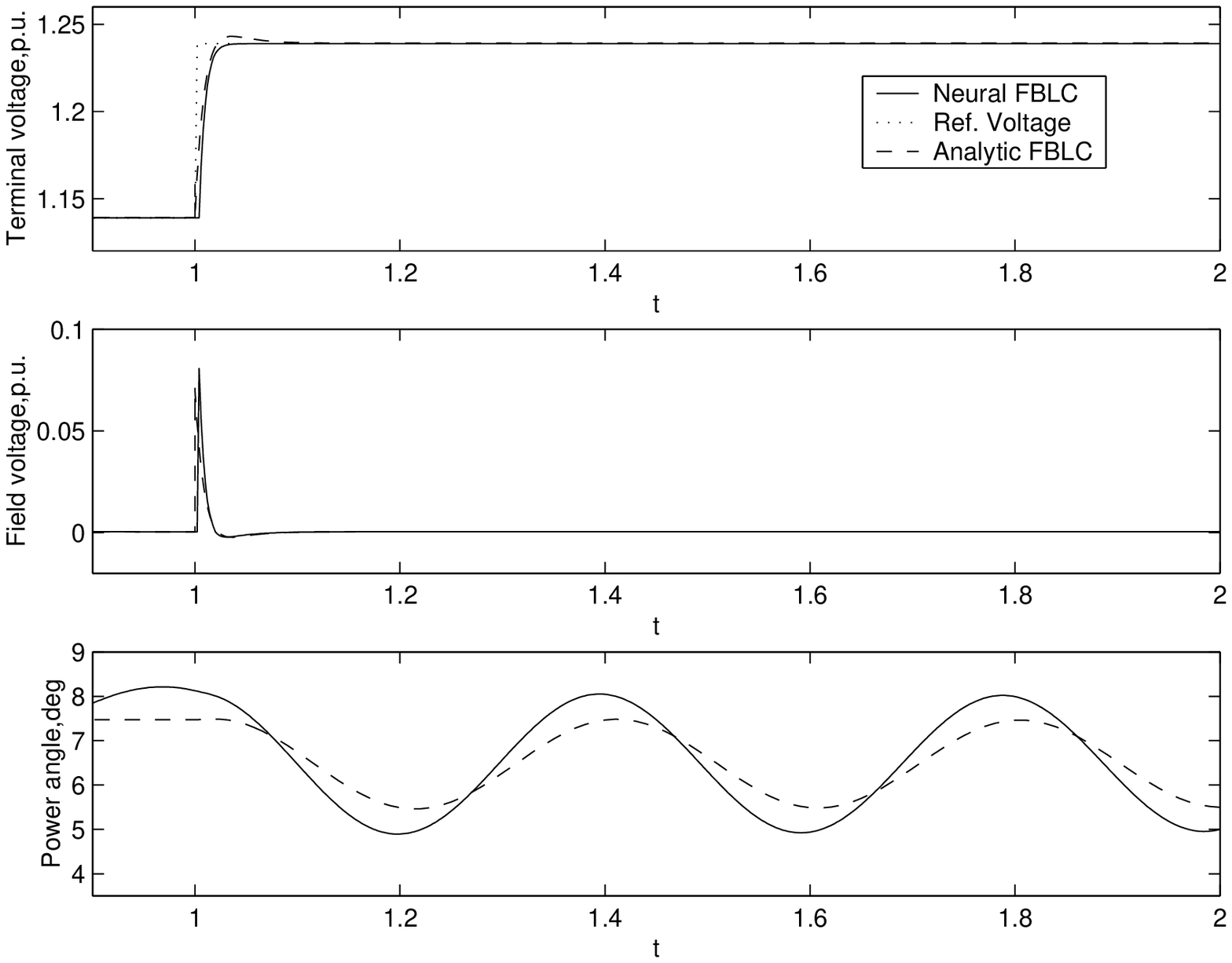,width=11cm}
\caption{Excitation Controller for Case $p=1$}
\label{case_p1}
\end{figure}
In order to do a fair comparison with the analytic feedback linearizing controller proposed in
\cite{dithesis}, we also compare the response for $p=3$ (since \cite{dithesis} used a simplified 3rd order
power system model for controller design). In this case, (to place the closed-loop poles at the same
locations as the analytic controller, $z \approx 0.7$), the
transfer function is of the form
\[ \frac{Y[z]}{R[z]} = \frac{K_{1}}{z^{3}+C_{2}z^{2}+C_{1}z+C_{0}}=\frac{0.03}{z^{3}-2.1z^{2}+1.47z-0.34} \]
The response is depicted in Fig. \ref{case_p3}. The neural controller is clearly better in this case since
the 3rd order analytical
feedback lineariztion did not decouple the high order modes.
\begin{figure}[h]
\centering\
\psfig{file=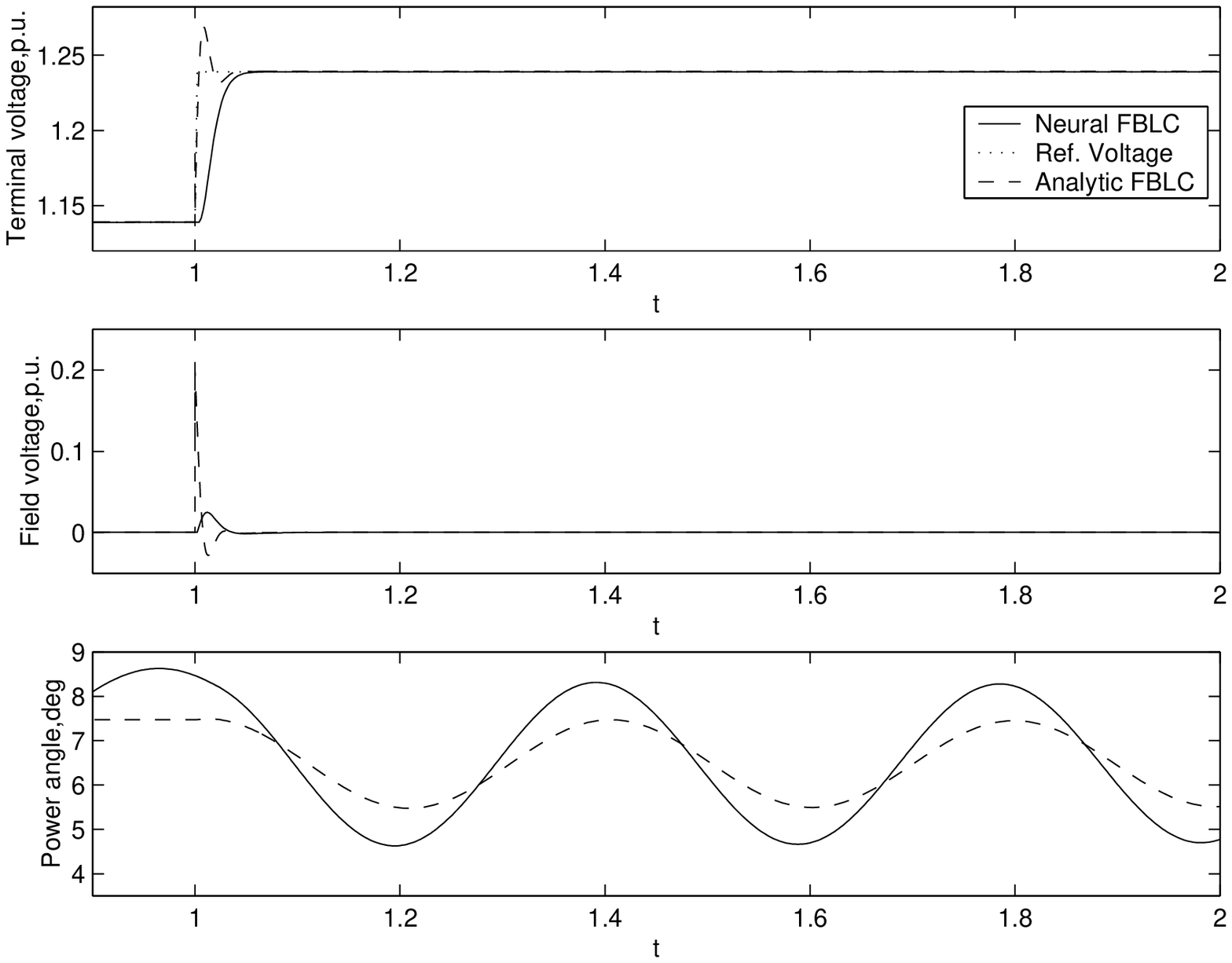,width=11cm}
\caption{Excitation Controller for Case $p=3$}
\label{case_p3}
\end{figure}

\item {\em Case, $p=7$}: \hspace{2mm}
the effects due to inexact pole/zero cancellation mean that we still have a $7^{th}$ order system to deal
with. It is proposed that a
controller be designed to ``force'' the system to behave like a stable 7th order system.
The control setup is shown in Fig. \ref{fig_kof1}.
\begin{figure}[h]
\centering\
\psfig{file=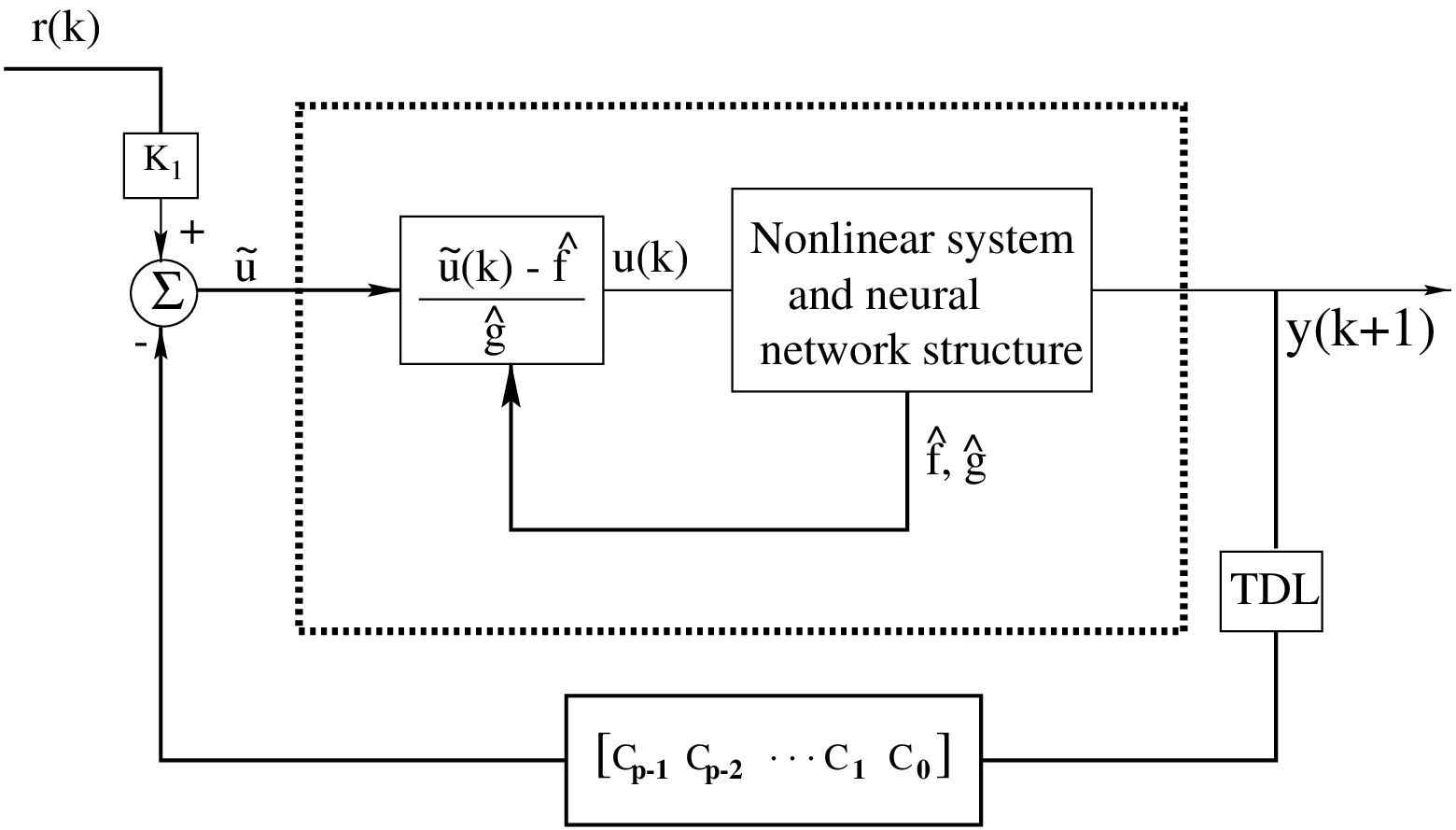,width=8cm}
\caption{Excitation Control setup for proposed approach (Case $p=7$)}
\label{fig_kof1}
\end{figure}
In practice, we know that exact pole/zero cancellation will not result from feedback linearization,
hence we will select p=7 and attempt to impose the desired $7^{th}$ order system behaviour and suppress 
the effects of inexact cancellation.
Let us choose to place the closed loop poles at $z=0.7$, to mimic the behaviour of the ST1A controller at 
the operating point. The resulting characteristic polynomial is given by
\begin{eqnarray*}
Q(z) & = & (z-0.7)^{7}  \\
         & = &  z^{7}-4.9z^{6}+10.29z^{5}-12.005z^{4}+8.4035z^{3}-3.5295z^{2}+0.8235z\\
         &   &  -0.0824
\end{eqnarray*}
Simulation results shown so far suggest that this later scheme gives better performance than the others
considered. The main challenge at this time is how to extend the same type of stability analysis and
convergence proofs
given in \cite{chen2} to this particular framework or develop alternative methods of analysis.
\end{itemize}
Now we compare the response of the proposed neurocontroller (with
$p=7$) 
the analytic
controller of \cite{dithesis} and the ST1A conventional controller.

\begin{figure}[h]
\centering\
\psfig{file=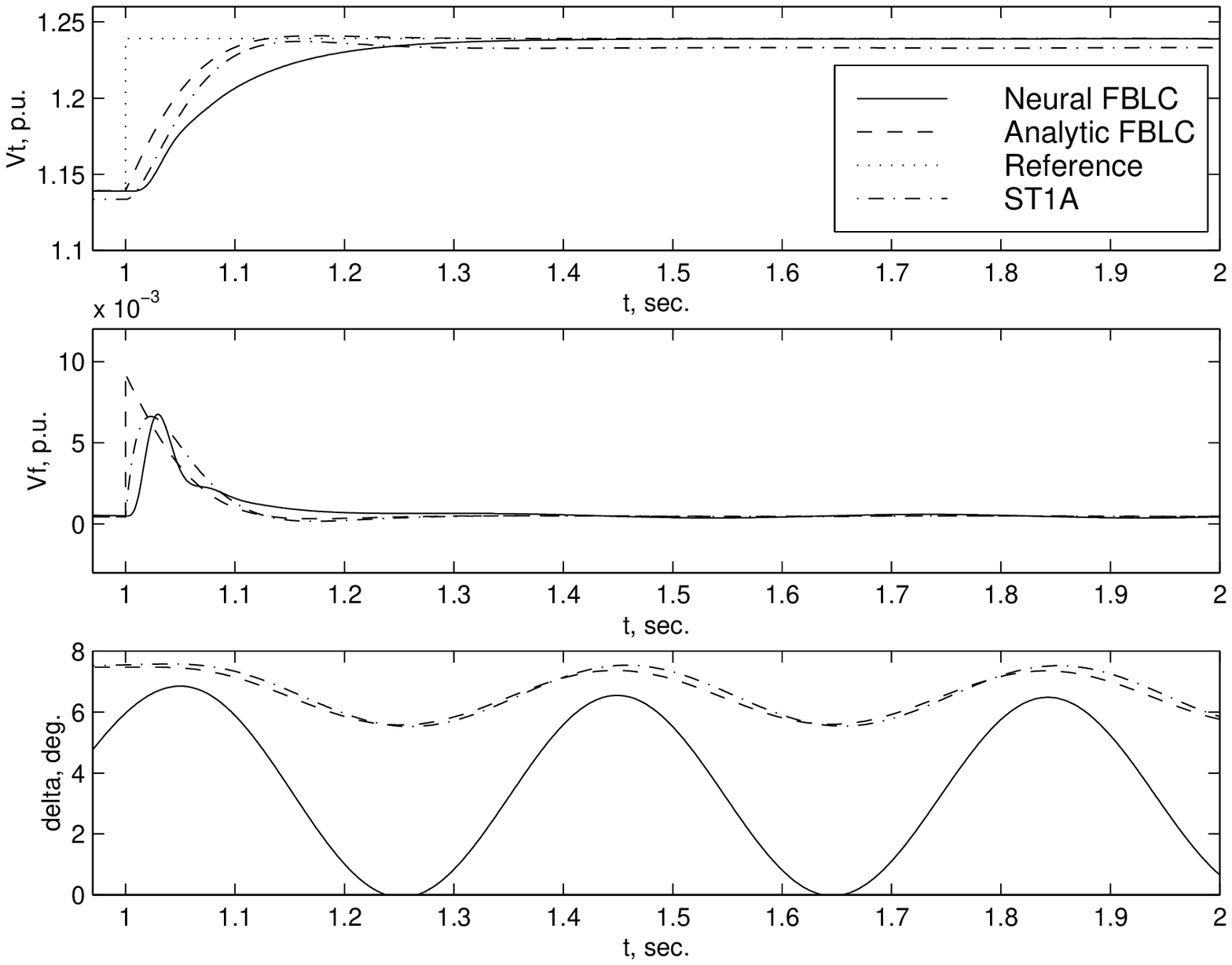,width=11cm}
\caption{Excitation control-perturbation from nominal operating point}
\label{all_excitation}
\end{figure}

The system output is initially at the nominal value of $v_{t} = 1.1392$. Then a step change of  $0.1\; p.u.$
in reference voltage occurs at $t = 1.0$ second. The
performance of all three controllers is depicted in Fig. \ref{all_excitation}. 
At $t=1.3$ seconds, $\frac{v_{t}}{v_{ref}}=0.9950$ for the ST1A exciter, $1.0002$ for the analytic controller and $0.9979$
for the neural controller. At $t=2.0$ seconds, $\frac{v_{t}}{v_{ref}}=0.9952$ for the ST1A exciter, $1.0000$ for the analytic
controller and $0.9998$ for the neural controller. The analytic controller initially gives a slight overshoot before settling
to the final value. The neural controller has a steady state error of $0.02 \% $ while the ST1A exciter gives an error of $0.5
\%$. 
The control effort is just less for the neural and ST1A excitation control than the analytic controller. 

\begin{figure}[h]
\centering\
\psfig{file=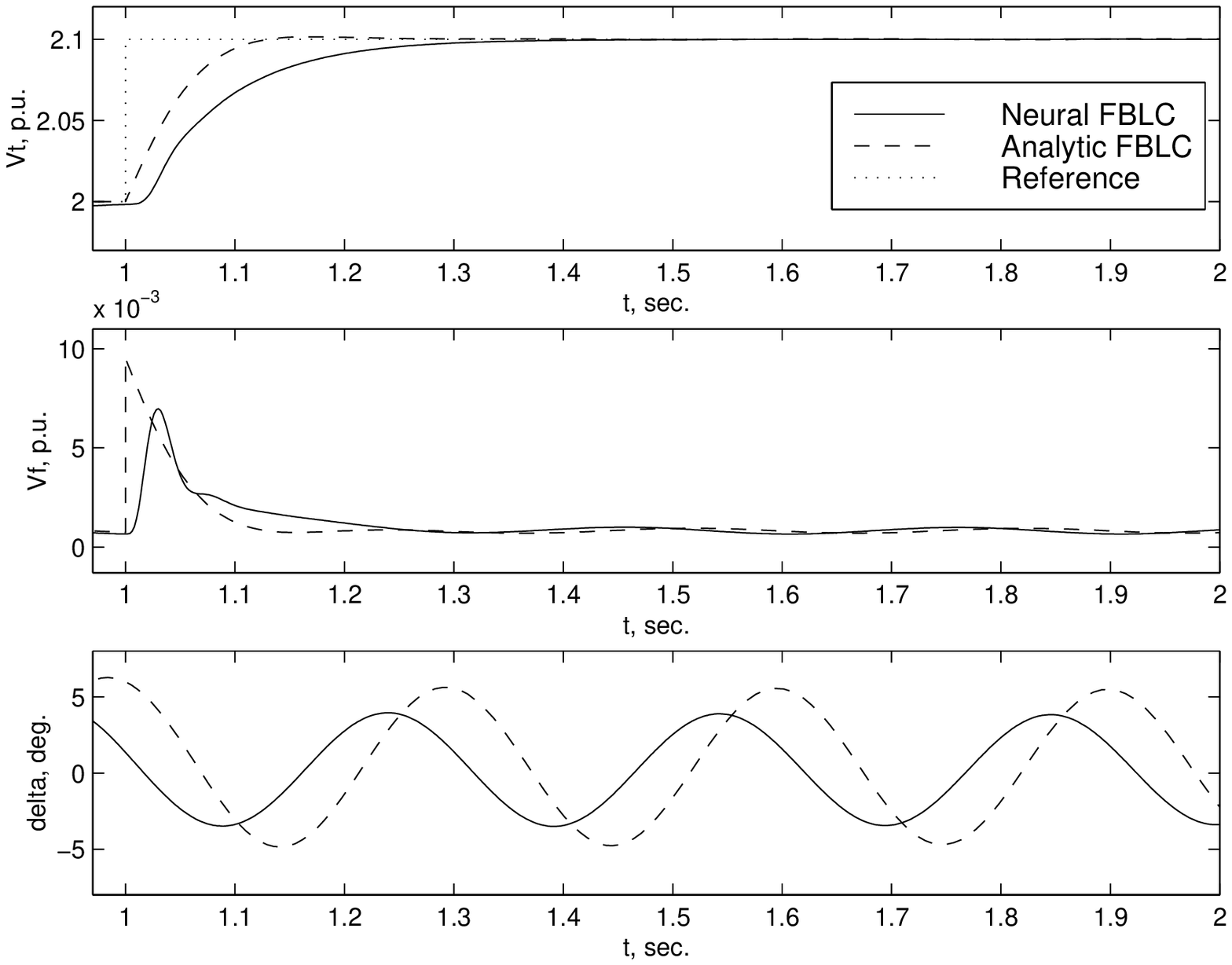,width=11cm}
\caption{Excitation control for change in operating point {\em far} from nominal}
\label{diff_excitation}
\end{figure}
Next we test the ability of  the neural and analytic excitation controllers to track reference signals at an
operating point far from the nominal, i.e. from $v_{ref}=2.0$ to $v_{ref}=2.1$. The step change takes place at $t  = 1.0$ second and the
response is depicted in Fig. \ref{diff_excitation}. It is observed that similar transient response
characteristics are
obtained compared to the response to a perturbation about the nominal operating point as shown in
Fig. \ref{all_excitation}.
Also the control effort $v_{f}$ is less for the neural controller than for the analytic controller. Finally, observe that in both
Fig. \ref{all_excitation} and Fig. \ref{diff_excitation}, the magnitude of the rotor angle oscillations
for the neural excitation
controller mimics the behaviour of the system under no control in response to a small perturbation. This is due to fact that feedback
linearization was applied to a $7^{th}$ order system, as opposed to a $3^{rd}$ order, as done with analytic feedback linearization.
Hence, the effect of excitation control is almost completely decoupled from the mechanical modes of the system.

{\bf Excitation controller/power system stabilizer:}
In \cite{dithesis}, the control law used for the power system stabilizer is given by
\[ v \; = \; -K(v_{t} - v_{ref} + K_{1} \dot{\delta}), \hspace{10mm} K, K_{1} > 0 \]
giving voltage tracking dynamics:
\[ \dot{v}_{t} \; = \; -K(v_{t} - v_{ref}) - K K_{1} \dot{\delta} \]
with gains
\[ K \; = \; 200 , \hspace{10mm} K_{1} \; = \; 0.7091 \]
These values are determined by a root locus criterion (see
\cite{dithesis}) in order to obtain optimal damping of rotor
angle oscillations.

\begin{figure}[h]
\centering\
\psfig{file=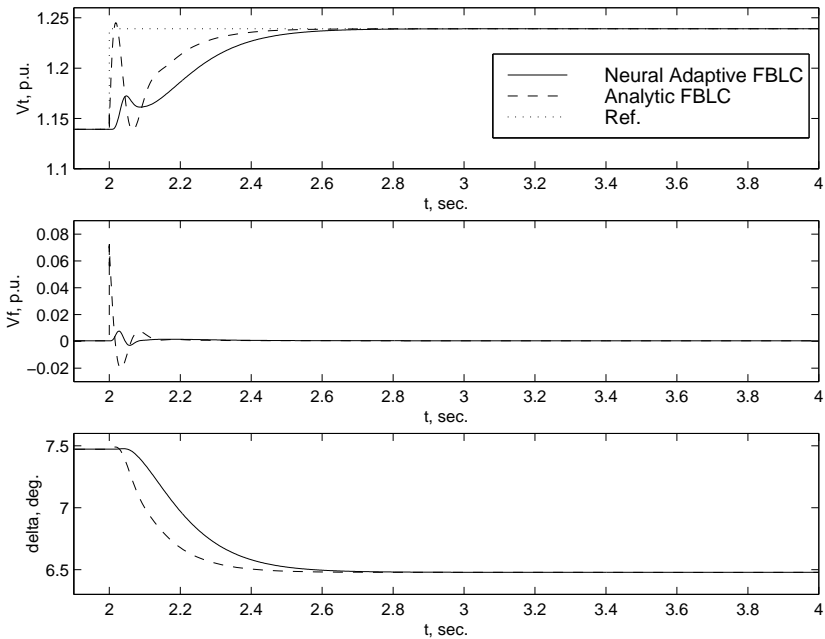,width=11cm} 
\caption{$v_{t},\;v_{f}\;$ and $\delta$ behaviour for a 0.1 p.u. step change in ref. voltage occuring at t=2.0
seconds} \label{perf_all}
\end{figure}

Recall that the control law for the neural adaptive controller/power system stabilizer is given by
\[ u_{k}=\frac{\tilde{u}(k)-\hat{f}(\cdot)}{\hat{g}(\cdot)} + K_{pss} \cdot \dot{\delta} (k) \]
where the gain $K_{pss}$ is given by
\begin{equation} K_{pss} \; = \; \nu \times 0.7091 . \label{nu_eqn}  \end{equation}
For the simulations which follow,
$ \nu \; = \; 3 .
$
In general,
increasing $\nu$ has the effect of improving the damping on
$\delta$ for step changes in reference signal. However $\nu$ can not be assigned arbitrarily large values. For instance, values of $\nu$
much
higher than $3$ have the effect of making the system unstable. 
The output of the system is originally at the nominal value of $v_{t}=1.1392$. A step change of $0.1 p.u.$ in
reference voltage occurs at $t=2.0$ seconds. Observe from Fig. \ref{perf_all} that the neural controller
accurately tracks the
reference signal (actually, tracking error is 0.0082 \% after 2.8 seconds) and quickly damps out the rotor angle oscillations.

\begin{figure}[h]
\centering\
\psfig{file=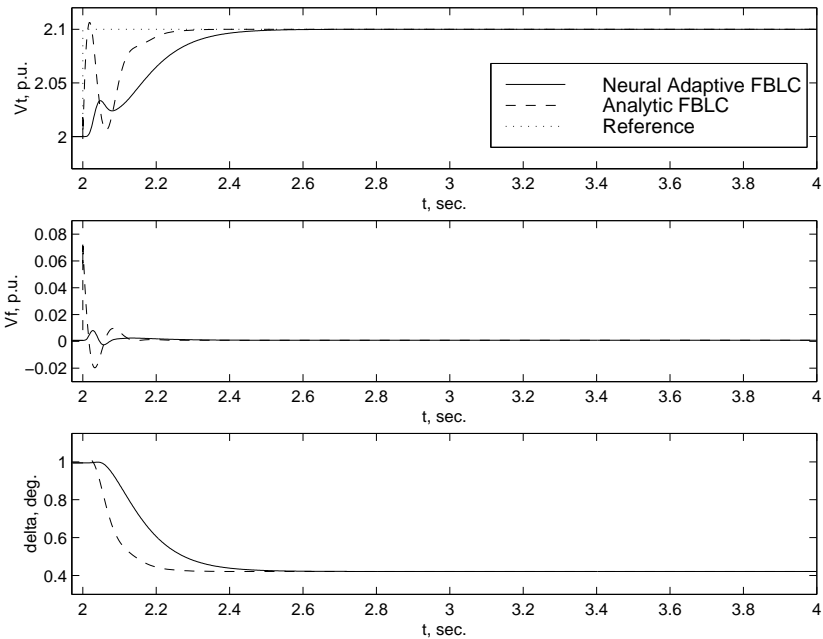,width=11cm}
\caption{$v_{t},\;v_{f}\;$ and $\delta$ behaviour for a 0.1 p.u. step change in ref. voltage occurring at t=2.0
seconds for operating point far from nominal}
\label{perf_all_diff}
\end{figure}

As done for the excitation controller, the controller/pss is validated for a change in operating point to 
$v_{ref}=2.0$ and a step change of $0.1$ p.u. takes place
at $t=2.0$ seconds. Dynamic behaviour similar to that in Fig. \ref{perf_all} are obtained for both
controllers. The figure
depicting the tracking and stabilizing dynamics for this new operating point is shown in
Fig. \ref{perf_all_diff}

\begin{figure}[h]
\centering\
\psfig{file=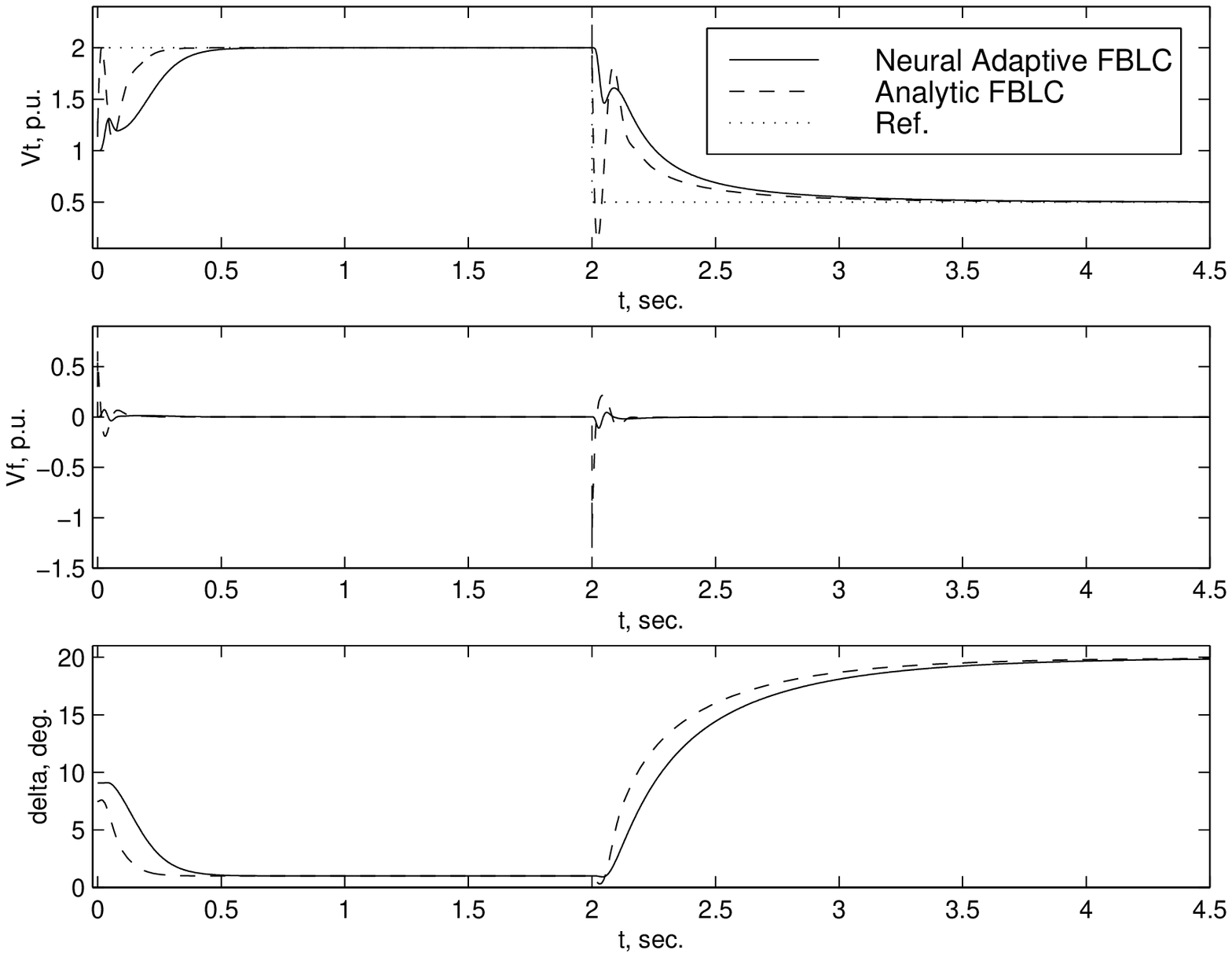,width=11cm}
\caption{$v_{t},\;v_{f}\;$ and $\delta$ behaviour for big changes in operating point}
\label{perf_all_plus}
\end{figure}
The neural controller was also tested for big changes in the operating points of the system as follows: 
a step increase in $v_{ref}$ from $1 \; p.u.$ to $2 \; p.u.$ was introduced at $t=0$. After $2$ seconds, $v_{ref}$ 
dropped to $0.5 \; p.u.$ Fig. \ref{perf_all_plus} depicts the results for both the neural and analytic
controllers. 
Both controllers performed well and exhibited similar tracking dynamics compared to the response to a small step change 
in reference voltage shown in Fig. \ref{perf_all}.  

The tuned analytic controller described above 
exhibits high frequency
harmonics during tracking due to the effects of the unmodelled dynamics in the design phase. The neural model does 
not exhibit such pronounced oscillations. Also observe that 
the control effort ($v_{f}$ in Figs. \ref{perf_all} and \ref{perf_all_plus}) is several orders of magnitude higher for the 
analytic controller than for the neural 
controller. Finally, observe that the steady state tracking error is effectively zero in both cases and the speed of 
response is comparable.
\begin{figure}[h]
\centering\
\psfig{file=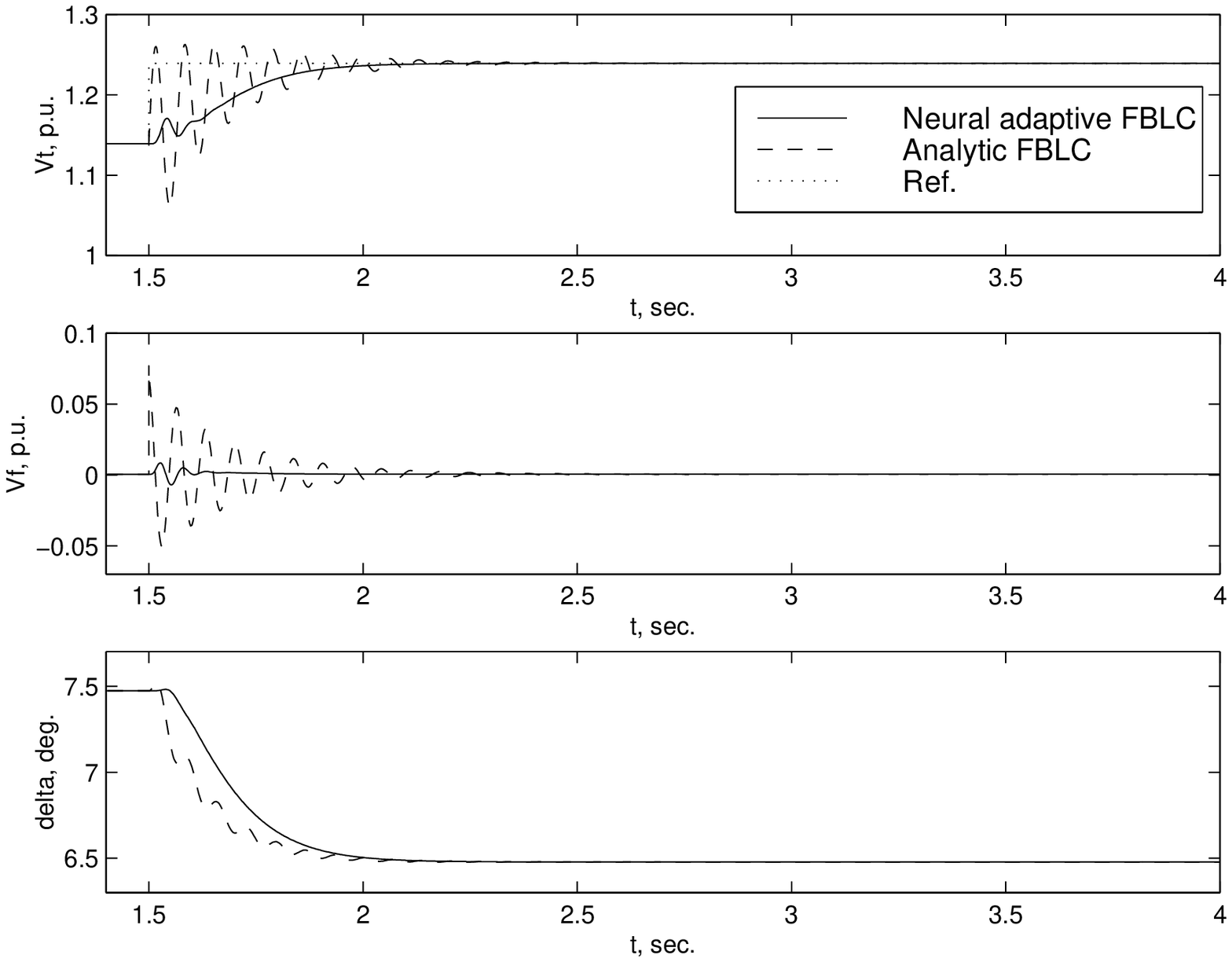,width=10cm}
\caption{$v_{t},\;v_{f}\;$ and $\delta$ behaviour under parameter variation}
\label{figs_vary}
\end{figure}
\begin{figure}[h]
\centering\
\psfig{file=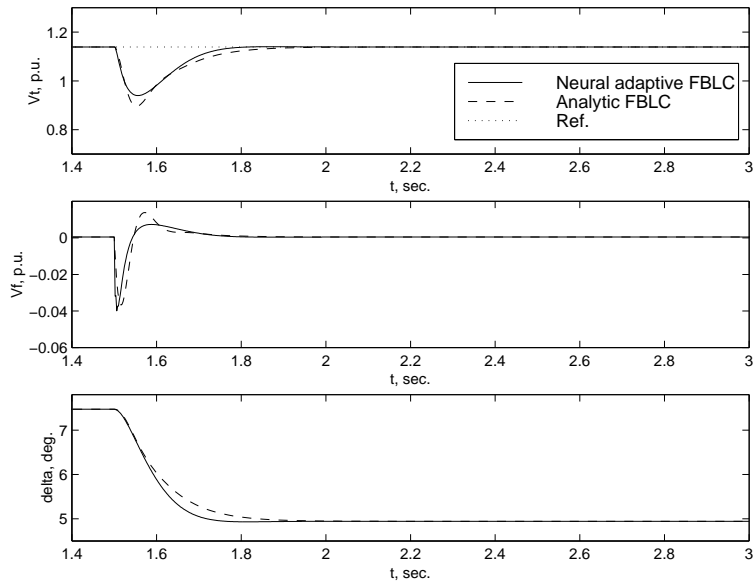,width=10cm}
\caption{$v_{t},\;v_{f}\;$ and $\delta$ behaviour under a disturbance}
\label{figs_disturb}
\end{figure}
\begin{figure}[h]
\centering\  
\psfig{file=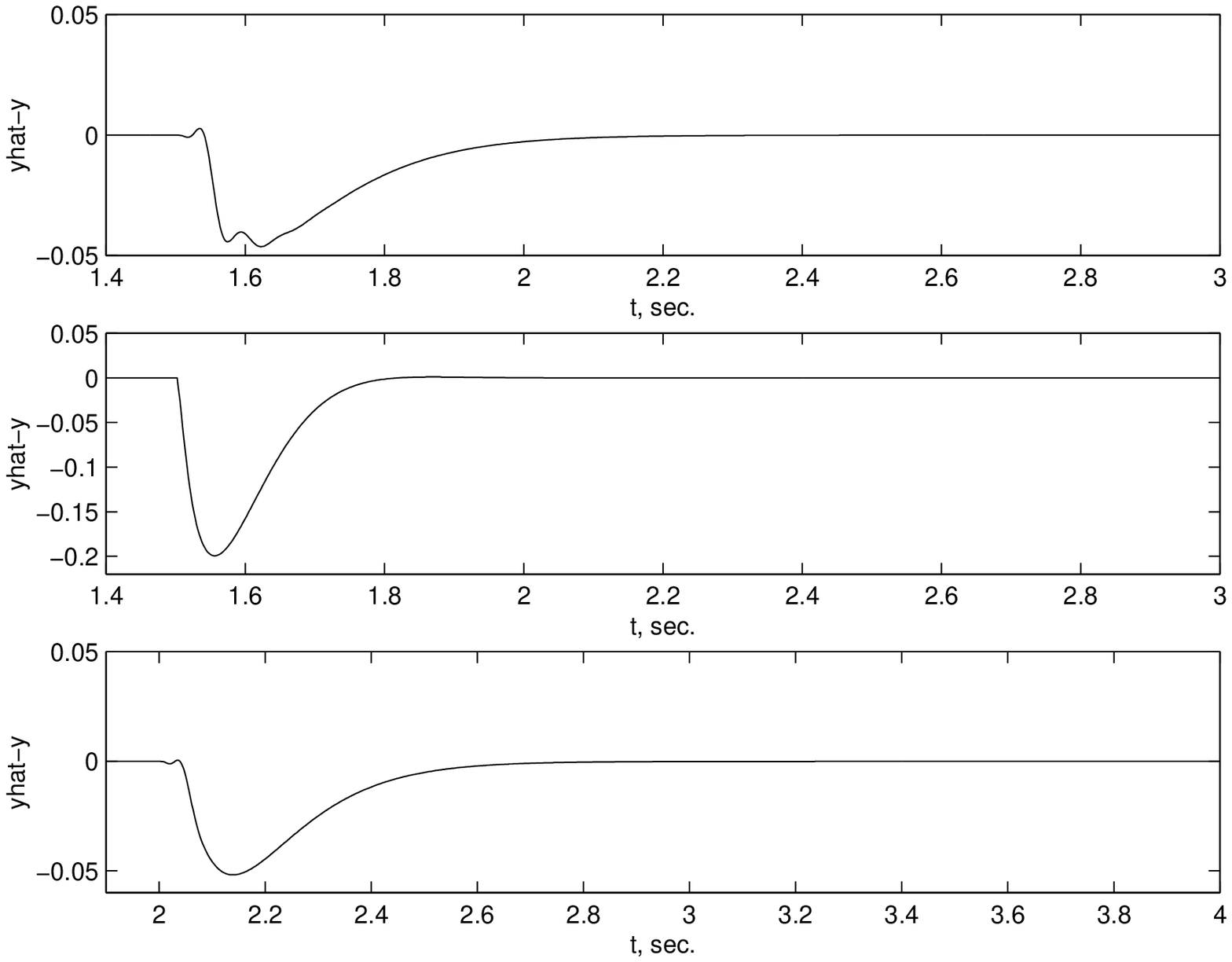,width=10cm}
\caption{Plant output errors for varying parameter(top), disturbance (middle) and excitation controller/pss
(bottom)}
\label{figs_with_error}
\end{figure}  
Another set of tests done on the control system involves the inertia parameter $H$  drifting from its nominal value of $9.5$
to about $4.75$ (i.e. a $50 \%$ drift) to see how well the control system behaves when feedback
linearization is inexact. This is depicted in Fig. \ref{figs_vary}.  The neural controller performs better
than the analytic controller in this case since the analytic controller design was based on
exact knowledge of all the parameters of the plant. It is observed that the tracking dynamics  of the neural
controller is essentially the same as what we have in Fig. \ref{perf_all}, where the parameter $H=9.5$.
For this test, the
value of $\nu$ in (\ref{nu_eqn}) is reduced to $\nu=2.6$. For arbitrary variation in the parameters, the neural controller
performance can be improved by
reducing the size of the deadzone used in the adaptation algorithm so that weights are updated for values of $d_{0}$ less than $0.01$
in this particular case.

The system performance under a disturbance
is
investigated by
forcing the mechanical power input,  $P_{m}$  to drop suddenly by $33.4 \%$
(from 1.6512 to 1.1); the response is shown in Fig. \ref{figs_disturb}. Observe that both controllers
recover quickly
and
restore the terminal voltage $v_{t}$ to its nominally value. The rotor angle changes to a new value but still retains its damping. In
this case both controllers gave the same level of performance, although the tracking dynamics for the neural controller is just
slightly better. We note that in practice this may be a rather drastic change in power level, but the
simulation is done just to see how far the system may be ``pushed'' whilst still maintaining control
action.

Finally, the effect of on-line weight adaptation on the response of the overall control system is investigated by observing the rate at
which the
error adaptation occurs. The Fig. \ref{figs_with_error} shows how quickly the plant output error
$\{\hat{y}(k+1)-y(k+1)\}$ is adapted
for
the
parameter variation,
disturbance and normal excitation controller/power system stabilizer tests. Observe that the errors are not too large
(apart from that associated with system disturbance) and are adapted quickly. Therefore it is concluded that adapting the weights online
does not unduly slow down the response of the overall control system.

\section{Concluding Remarks} 
A neural adaptive feedback linearizing excitation controller/ power system stabilizer  has been proposed for 
a high-order model of the synchronous machine/infinite bus power system. The controller provides 
tracking performance with less high frequency harmonics than an analytic controller based on a reduced-order model of the system
and would seem to tolerate parameter drift. In 
the cases considered, 
step changes in reference voltage are tracked, with the neural controller requiring significantly less control 
effort than the analytic controller. The neural controller also provides damping of power angle 
oscillations with similar voltage and power angle transient characteristics over different operating points of the 
system. Also on-line controller adaptation does not slow down the overall control system response since output errors are reduced
quickly. The advantage of this approach lies in the fact
that exact knowledge of the power system dynamics is not absolutely
necessary
and the system state need not all be measurable. Also, the complicated
mathematical mapping required to carry out
analytic feedback linearization is avoided.
The control scheme is quite straightforward and easy to synthesize. Thorough mathematical analyses
guaranteeing convergence and/or stability for this type of control scheme (comparable to that
reported in \cite{chen2}) is still to be done.

\end{document}